# SEMI-ANALYTICAL MASS MATRIX FOR 8-NODE BRICK ELEMENT.

E Hanukah
Faculty of Mechanical Engineering, Technion – Israel Institute of Technology, Haifa 32000, Israel
Corresponding author Email: eliezerh@tx.technion.ac.il

**Abstract**

Nowadays integration of mass matrix components in the element domain is performed using various numerical integration schemes, each one possess different level of accuracy, alters in number of integration (Gauss) points and requires different amount of computations. Herein semi-analytical approach is suggested. Metric (Jacobian determinant) is approximately modeled using its evaluations in certain points. Analytical integration is performed to derive simple explicit closed-form expressions for each term of the mass matrix. Two schemes are discussed: the first assumes constant metric (CM) in the initial domain, using evaluation at the centroid. The second allows for linear variation of the metric (LM linear metric) in the domain using 3 additional evaluation points. Both schemes are exact for rectangular and non-rectangular parallelepiped mesh. Careful symbolic manipulations and convenient choice of evaluation points allow us to avoid unnecessary operations. The accuracy of both schemes is studied numerically using randomly generated coarse mesh. Significant superiority in accuracy over equivalent schemes is reported. An important implication of this study is that it can replace currently used schemes.

**Key words**: hexahedral element, consistent / lumped mass matrix, closed-form, symbolic computational mechanics.

## 1. Introduction

Probably every book or lecture notes concerning finite element method (FEM) for solid continuum includes 8 node brick element, as well as roughly all commercial widely-used packages e.g. ABAQUS$^{TM}$, ANSYS$^{TM}$, LS-DYNA$^{TM}$ etc., has it implemented. Therefore, derivation of sufficiently accurate and computationally inexpensive integration rule for consistent and lumped mass matrix is vital.

Mass matrix components, internal forces and stiffness matrix, all require integration in the element domain, which is obtained with the help on numerical integration schemes e.g. [1-4]. Several studies exist that exploit the idea of analytical integration for stiffness matrixes resulting in greater accuracy and efficiency e.g. [5-8]. Furthermore, hierarchical semi-analytical displacement based approach is used to model three dimensional bodies e.g. [9-12] yielding analytical and numerical solutions.

In present study semi-analytical approach for computation of mass matrix components is offered. Approximation for the metric (Jacobian determinant) is formulated based on its evaluation at certain points and polynomial coordinates dependence. Analytical integration is performed to derive explicit closed-form expressions for the mass matrix components in terms of initial densities and metric evaluations.

Our first CM rule assumes constant metric in the domain which is sampled at the centroid. Closed-form expressions for mass matrix components follow from analytical integration. Our second LM rule assumes linear variation of the metric in the domain. To this end three additional evaluation points are required. We show that the exact metric is fourth order with respect to coordinates; hence, the considered CM and LM are low order schemes.

Preliminary numerical study is conducted to test the performance of new rules. CM and LM are exact for parallelepiped mesh. Random coarse mesh elements are generated and the averaged absolute error is calculated with respect to exact results. It is found that CM is significantly over performs numerical integration based on one point quadrature. While LM is superior to numerical integration scheme based on 4 integration points. The study considers commonly used lumped matrix formulation; however the extension to consistent mass matrix is straightforward.

The outline of the paper is as follows. Section 2 recalls necessary definitions and basics of 8-node brick element such as the shape functions, kinematic approximation, initial density approximation, mass matrix definition etc. Section 3 presents all the details of the proposed integration rules applied to widely-used lumped mass matrix formulation. Section 4 contains preliminary numerical accuracy study, including comparison to equivalent schemes. Section 5 records our conclusions.

## 2. Background.

Initial location of the nodes of the standard 8-node brick element (e.g. [13] pp.68) is denoted by $\mathbf{N}_i\,(i=1,..,8)$, its components are given in terms of global Cartesian coordinates system $\mathbf{N}_i = N_{ik}\mathbf{e}_k\,(i=1,..,8, k=1,2,3)$, traditional summation convention on repeated index is implied. The shape functions $\varphi^i\,(i=1,..,8)$ in terms of local convected coordinate system $\{\xi,\eta,\zeta\}$ is given by

$$
\begin{aligned}
\varphi^1 &= (1-\zeta-\eta+\eta\zeta-\xi+\xi\zeta+\xi\eta-\xi\eta\zeta)/8\\
\varphi^2 &= (1-\zeta-\eta+\eta\zeta+\xi-\xi\zeta-\xi\eta+\xi\eta\zeta)/8\\
\varphi^3 &= (1-\zeta+\eta-\eta\zeta+\xi-\xi\zeta+\xi\eta-\xi\eta\zeta)/8\\
\varphi^4 &= (1-\zeta+\eta-\eta\zeta-\xi+\xi\zeta-\xi\eta+\xi\eta\zeta)/8\\
\varphi^5 &= (1+\zeta-\eta-\eta\zeta-\xi-\xi\zeta+\xi\eta+\xi\eta\zeta)/8\\
\varphi^6 &= (1+\zeta-\eta-\eta\zeta+\xi+\xi\zeta-\xi\eta-\xi\eta\zeta)/8\\
\varphi^7 &= (1+\zeta+\eta+\eta\zeta+\xi+\xi\zeta+\xi\eta+\xi\eta\zeta)/8\\
\varphi^8 &= (1+\zeta+\eta+\eta\zeta-\xi-\xi\zeta-\xi\eta-\xi\eta\zeta)/8
\end{aligned} \tag{1}
$$

Material point $X$ occupies location $\mathbf{X}$ inside the element domain $-1 \le \xi, \eta, \zeta \le +1$ is given by

$$\mathbf{X} = \varphi^i \mathbf{N}_i \ \ (i = 1,..,8) \tag{2}$$

The initial density of material points initially located at the nodes are denoted by $\rho_i \, (i = 1,..,8)$ and the density inside the domain is approximated by

$$
\begin{aligned}
&\rho_0 = \varphi^i \rho_i \ \ (i = 1,..,8)\\
&\sum_{i=1}^{8} \varphi^i = 1 \ , \ \rho_0\big|_{\rho_1=\rho_2=\rho_3=\ldots=\rho_8=\rho^*} = \rho^*
\end{aligned} \tag{3}
$$

Where $\rho_0(\xi, \eta, \zeta)$ stands for the initial density, the above approximation admits homogeneity in the case of constant initial density at nodes.

The jacobian determinant (metric) of global-local coordinates transformation $J$, differential volume element $dV$, and initial volume $V$ are defined by

$$
\begin{aligned}
&J = \mathbf{X}_{,1} \times \mathbf{X}_{,2} \bullet \mathbf{X}_{,3} = \begin{vmatrix} (\mathbf{X}\bullet\mathbf{e}_1)_{,1} & (\mathbf{X}\bullet\mathbf{e}_1)_{,2} & (\mathbf{X}\bullet\mathbf{e}_1)_{,3} \\ (\mathbf{X}\bullet\mathbf{e}_2)_{,1} & (\mathbf{X}\bullet\mathbf{e}_2)_{,2} & (\mathbf{X}\bullet\mathbf{e}_2)_{,3} \\ (\mathbf{X}\bullet\mathbf{e}_3)_{,1} & (\mathbf{X}\bullet\mathbf{e}_3)_{,2} & (\mathbf{X}\bullet\mathbf{e}_3)_{,3} \end{vmatrix} > 0 \ , \ J_{mn}(\xi, \eta, \zeta, N_{ik}) = (\mathbf{X}\bullet\mathbf{e}_m)_{,n}\\
&dV = Jd\xi d\eta d\zeta = Jd\square \ , \ V = \int_{-1}^{+1}\int_{-1}^{+1}\int_{-1}^{+1} dV \ , \ (i = 1,..,8, m, n, k = 1, 2, 3)
\end{aligned} \tag{4}
$$

Where $(\times)$ and $(\bullet)$ stand for vector cross and scalar products and $|\bullet|$ stand for determinant operator, comma denotes partial differentiation with respect to coordinates. Here and throughout the study, determinant of general (non-symmetric) 3x3 matrix is computed as

$$J := J_{11}J_{22}J_{33} - J_{11}J_{23}J_{32} - J_{31}J_{22}J_{13} - J_{21}J_{12}J_{33} + J_{21}J_{32}J_{13} + J_{31}J_{12}J_{23} \tag{5}$$

The above is consistent with standard definition of determinant. Isoparametric formulation (e.g.[14] pp.104) for mass conserving element, yield the next consistent, symmetric, and positive definite mass matrix

$$M^{ij} = \int_V \rho_0 \phi^i \phi^j dV = \int_{-1}^{+1}\int_{-1}^{+1}\int_{-1}^{+1} \rho_0(\xi,\eta,\zeta,\rho_r)\phi^i(\xi,\eta,\zeta)\phi^j(\xi,\eta,\zeta)J(\xi,\eta,\zeta,N_{sk})d\xi d\eta d\zeta \tag{6}$$

$$M^{ij} = M^{ji} \ , \ (i,j,r,s = 1,..,8, k = 1,2,3)$$

Lumped (diagonal) mass matrix is preferred in explicit (and mostly for implicit) transient analysis using 8-noded brick (e.g. [13] pp.140). Several forms have been suggested, however, the commonly used one (e.g. [15, 16]) is given by

$$M^{ij} = 0 \ \ (i \neq j) \ , \ M^{ii} = \int_V \rho_0 \phi^i dV \ , \ (i,j = 1,..,8) \tag{7}$$

Due to the above practical reason, herein we focus on the lumped formulation, nevertheless all the same procedures applicable as well for a consistent mass matrix (6).

With the help of Taylor's multivariable expansion about the centroid of an element $\mathbf{X}_0 = \mathbf{X}(\xi = 0, \eta = 0, \zeta = 0)$, one can exactly represent $J$ as

$$\begin{aligned} J = {} & J_0 + \\ & \xi J_1 + \eta J_2 + \zeta J_3 + \\ & \xi\eta J_4 + \xi\zeta J_5 + \eta\zeta J_6 + \xi^2 J_7 + \eta^2 J_8 + \zeta^3 J_9 + \\ & \xi\eta\zeta J_{10} + \xi^2\eta J_{11} + \xi\eta^2 J_{12} + \xi^2\zeta J_{13} + \eta^2\zeta J_{14} + \xi\zeta^2 J_{15} + \eta\zeta^2 J_{16} + \\ & \xi^2\eta\zeta J_{17} + \xi\eta^2\zeta J_{18} + \xi\eta\zeta^2 J_{19} \end{aligned} \tag{8}$$

$$J_0 = J\big|_{\xi,\eta,\zeta=0} = J(\mathbf{X}_0) \ , \ J_1 = \frac{\partial J(\mathbf{X}_0)}{\partial \xi} \ , \ J_2 = \frac{\partial J(\mathbf{X}_0)}{\partial \eta} \ , \ J_3 = \frac{\partial J(\mathbf{X}_0)}{\partial \zeta} \ ,..., \ J_{19} = \frac{\partial^4 J(\mathbf{X}_0)}{2\partial\zeta^2 \partial\xi \partial\eta}$$

The metric $J$ is fourth order with respect to coordinates. It is important to emphasize that for parallelepiped mesh, the metric is independent of coordinates $J = J_0$ (constant metric). Using the above representation (8) together with (1) and (3) analytical integration of the lumped mass matrix component (7) is performed and used later as an exact values with respect to which the error is computed. Here and throughout the study, computer algebra system (CAS) MAPLE™ were used to perform all the symbolic manipulations, including integration, differentiation, simplification, direct translation of explicit expression to Fortran77, generation of random numbers etc.

Standard numerical integration in element domain is recalled (e.g. [14] pp.121)

$$\int_{-1}^{+1}\int_{-1}^{+1}\int_{-1}^{+1} fJd\Box = \sum_{p=1}^{n_p} f(\xi_p,\eta_p,\zeta_p)J(\xi_p,\eta_p,\zeta_p)w_p \tag{9}$$

Where $n_p$ stand for number of integration points, $w_p$ denotes weights at integration points and $\xi_p, \eta_p, \zeta_p$ are coordinates of integration points. Special integration for hexahedral (brick) elements is given by

$$n_p = 1 \Rightarrow (\xi_1, \eta_1, \zeta_1, w_1) = (0,0,0,8)$$

$$n_4 = 4 \Rightarrow (\xi_1, \eta_1, \zeta_1, w_1) = (0, \sqrt{2/3}, -1/\sqrt{3}, 2) \;, \; (\xi_2, \eta_2, \zeta_2, w_2) = (0, -\sqrt{2/3}, -1/\sqrt{3}, 2) \quad (10)$$
$$(\xi_3, \eta_3, \zeta_3, w_3) = (\sqrt{2/3}, 0, 1/\sqrt{3}, 2) \;, \; (\xi_4, \eta_4, \zeta_4, w_4) = (-\sqrt{2/3}, 0, 1/\sqrt{3}, 2)$$

Where $n_p = 1$ is one point integration rule and $n_p = 4$ is four point integration rule.

## 3. Semi-analytical approach.

For the first CM rule we neglect all coordinate dependent terms in (8) and approximate the metric by constant $J \approx J_0$

$$J_0 = J(\xi{=}0, \eta = 0, \zeta{=}0) {=} \det(J^0_{mn}) \quad (11)$$

Components $J^0_{mn}$ are given by

$$\begin{aligned}
J^0_{11} &= 0.125(-N_{1,1}+N_{2,1}+N_{3,1}-N_{4,1}-N_{5,1}+N_{6,1}+N_{7,1}-N_{8,1})\\
J^0_{12} &= 0.125(-N_{1,1}-N_{2,1}+N_{3,1}+N_{4,1}-N_{5,1}-N_{6,1}+N_{7,1}+N_{8,1})\\
J^0_{13} &= 0.125(-N_{1,1}-N_{2,1}-N_{3,1}-N_{4,1}+N_{5,1}+N_{6,1}+N_{7,1}+N_{8,1})\\
J^0_{21} &= 0.125(-N_{1,2}+N_{2,2}+N_{3,2}-N_{4,2}-N_{5,2}+N_{6,2}+N_{7,2}-N_{8,2})\\
J^0_{22} &= 0.125(-N_{1,2}-N_{2,2}+N_{3,2}+N_{4,2}-N_{5,2}-N_{6,2}+N_{7,2}+N_{8,2}) \qquad (12)\\
J^0_{23} &= 0.125(-N_{1,2}-N_{2,2}-N_{3,2}-N_{4,2}+N_{5,2}+N_{6,2}+N_{7,2}+N_{8,2})\\
J^0_{31} &= 0.125(-N_{1,3}+N_{2,3}+N_{3,3}-N_{4,3}-N_{5,3}+N_{6,3}+N_{7,3}-N_{8,3})\\
J^0_{32} &= 0.125(-N_{1,3}-N_{2,3}+N_{3,3}+N_{4,3}-N_{5,3}-N_{6,3}+N_{7,3}+N_{8,3})\\
J^0_{32} &= 0.125(-N_{1,3}-N_{2,3}-N_{3,3}-N_{4,3}+N_{5,3}+N_{6,3}+N_{7,3}+N_{8,3})
\end{aligned}$$

Analytical integration is performed resulting in the next CM lumped mass matrix

$$
\begin{aligned}
M^{11} &= (8\rho_1 + 4\rho_2 + 2\rho_3 + 4\rho_4 + 4\rho_5 + 2\rho_6 + \rho_7 + 2\rho_8)J_0 / 27 \\
M^{22} &= (4\rho_1 + 8\rho_2 + 4\rho_3 + 2\rho_4 + 2\rho_5 + 4\rho_6 + 2\rho_7 + \rho_8)J_0 / 27 \\
M^{33} &= (2\rho_1 + 4\rho_2 + 8\rho_3 + 4\rho_4 + \rho_5 + 2\rho_6 + 4\rho_7 + 2\rho_8)J_0 / 27 \\
M^{44} &= (4\rho_1 + 2\rho_2 + 4\rho_3 + 8\rho_4 + 2\rho_5 + \rho_6 + 2\rho_7 + 4\rho_8)J_0 / 27 \\
M^{55} &= (4\rho_1 + 2\rho_2 + \rho_3 + 2\rho_4 + 8\rho_5 + 4\rho_6 + 2\rho_7 + 4\rho_8)J_0 / 27 \\
M^{66} &= (2\rho_1 + 4\rho_2 + 2\rho_3 + \rho_4 + 4\rho_5 + 8\rho_6 + 4\rho_7 + 2\rho_8)J_0 / 27 \\
M^{77} &= (\rho_1 + 2\rho_2 + 4\rho_3 + 2\rho_4 + 2\rho_5 + 4\rho_6 + 8\rho_7 + 4\rho_8)J_0 / 27 \\
M^{88} &= (2\rho_1 + \rho_2 + 2\rho_3 + 4\rho_4 + 4\rho_5 + 2\rho_6 + 4\rho_7 + 8\rho_8)J_0 / 27
\end{aligned}
\tag{13}
$$

Next we allow linear variation of the metric. Exact partial derivatives $J_1, J_2, J_3$ given by (8) are computed as explicit functions of nodal positions $N_{im}$. However they are found to be rather lengthy, each includes 194 additive terms while each term include 3 multiplications. Consequently approximation is suggested

$$J \approx J_0 + \xi\tilde{J}_1 + \eta\tilde{J}_2 + \zeta\tilde{J}_3 \tag{14}$$

The additional (first order) terms $\tilde{J}_k \ (k = 1,2,3)$ is given by

$$\tilde{J}_k = J_{point\,k} - J_0 \ , \ (k = 1,2,3) \tag{15}$$

Where $J_{point\,k}$ are the metric evaluations at 3 convenient points, which keep the components short

$$
\begin{aligned}
J_{point\,1} &= J(\xi{=}1, \eta = 0, \zeta{=}0) = \det(J^1_{mn}) \\
J_{point\,2} &= J(\xi{=}0, \eta = 1, \zeta{=}0) = \det(J^2_{mn}) \\
J_{point\,3} &= J(\xi{=}0, \eta = 0, \zeta{=}1) = \det(J^3_{mn}) \ \ (m, n = 1,2,3)
\end{aligned}
\tag{16}
$$

Components $J^k_{mn} \ (k, m, n = 1,2,3)$ are given by

$$
\begin{aligned}
& J^1_{11} = 0.125(-N_{1,1} + N_{2,1} + N_{3,1} - N_{4,1} - N_{5,1} + N_{6,1} + N_{7,1} - N_{8,1}) \\
& J^1_{12} = 0.25(-N_{2,1} + N_{3,1} - N_{6,1} + N_{7,1}) \ , \ J^1_{13} = 0.25(-N_{2,1} - N_{3,1} + N_{6,1} + N_{7,1}) \\
& J^1_{21} = 0.125(-N_{1,2} + N_{2,2} + N_{3,2} - N_{4,2} - N_{5,2} + N_{6,2} + N_{7,2} - N_{8,2}) \\
& J^1_{22} = 0.25(-N_{2,2} + N_{3,2} - N_{6,2} + N_{7,2}) \ , \ J^1_{23} = 0.25(-N_{2,2} - N_{3,2} + N_{6,2} + N_{7,2}) \\
& J^1_{31} = 0.125(-N_{1,3} + N_{2,3} + N_{3,3} - N_{4,3} - N_{5,3} + N_{6,3} + N_{7,3} - N_{8,3}) \\
& J^1_{32} = 0.25(-N_{2,3} + N_{3,3} - N_{6,3} + N_{7,3}) \ , \ J^1_{32} = 0.25(-N_{2,3} - N_{3,3} + N_{6,3} + N_{7,3})
\end{aligned}
\tag{17}
$$

$$J^2_{11}=0.25(N_{3,1}-N_{4,1}+N_{7,1}-N_{8,1})\ ,\ J^2_{13}=0.25(-N_{3,1}-N_{4,1}+N_{7,1}+N_{8,1})$$
$$J^2_{12}=0.125(-N_{1,1}-N_{2,1}+N_{3,1}+N_{4,1}-N_{5,1}-N_{6,1}+N_{7,1}+N_{8,1})$$
$$J^2_{21}=0.25(N_{3,2}-N_{4,2}+N_{7,2}-N_{8,2})\ ,\ J^2_{23}=0.25(-N_{3,2}-N_{4,2}+N_{7,2}+N_{8,2})$$
$$J^2_{22}=0.125(-N_{1,2}-N_{2,2}+N_{3,2}+N_{4,2}-N_{5,2}-N_{6,2}+N_{7,2}+N_{8,2}) \tag{18}$$
$$J^2_{31}=0.25(N_{3,3}-N_{4,3}+N_{7,3}-N_{8,3})\ ,\ J^2_{32}=0.25(-N_{3,3}-N_{4,3}+N_{7,3}+N_{8,3})$$
$$J^2_{32}=0.125(-N_{1,3}-N_{2,3}+N_{3,3}+N_{4,3}-N_{5,3}-N_{6,3}+N_{7,3}+N_{8,3})$$

$$J^3_{11}=0.25(-N_{5,1}+N_{6,1}+N_{7,1}-N_{8,1})\ ,\ J^3_{12}=0.25(-N_{5,1}-N_{6,1}+N_{7,1}+N_{8,1})$$
$$J^3_{13}=0.125(-N_{1,1}-N_{2,1}-N_{3,1}-N_{4,1}+N_{5,1}+N_{6,1}+N_{7,1}+N_{8,1})$$
$$J^3_{21}=0.25(-N_{5,2}+N_{6,2}+N_{7,2}-N_{8,2})\ ,\ J^3_{22}=0.25(-N_{5,2}-N_{6,2}+N_{7,2}+N_{8,2})$$
$$J^3_{23}=0.125(-N_{1,2}-N_{2,2}-N_{3,2}-N_{4,2}+N_{5,2}+N_{6,2}+N_{7,2}+N_{8,2}) \tag{19}$$
$$J^3_{31}=0.25(-N_{5,3}+N_{6,3}+N_{7,3}-N_{8,3})\ ,\ J^3_{32}=0.25(-N_{5,3}-N_{6,3}+N_{7,3}+N_{8,3})$$
$$J^3_{32}=0.125(-N_{1,3}-N_{2,3}-N_{3,3}-N_{4,3}+N_{5,3}+N_{6,3}+N_{7,3}+N_{8,3})$$

Using approximation (14) combined with analytical integration, the lumped mass matrix (7) turn out to be

$$\begin{aligned}
27M^{11} &= \left(-\tilde{J}_2-\tilde{J}_1-\tilde{J}_3\right)4\rho_1+\left(2\rho_8+\rho_7+2\rho_6+4\rho_5+4\rho_4+2\rho_3+4\rho_2+8\rho_1\right)J_0+\\
&\quad +\left(-\rho_8-2\rho_4-2\rho_5\right)\tilde{J}_1+\left(-2\rho_2-2\rho_5-\rho_6\right)\tilde{J}_2-\left(2\rho_2+\rho_3+2\rho_4\right)\tilde{J}_3\\
27M^{22} &= \left(-\tilde{J}_3-\tilde{J}_2+\tilde{J}_1\right)4\rho_2+\left(\rho_8+2\rho_7+4\rho_6+2\rho_5+2\rho_4+4\rho_3+8\rho_2+4\rho_1\right)J_0+\\
&\quad +\left(\rho_7+2\rho_3+2\rho_6\right)\tilde{J}_1+\left(-2\rho_1-\rho_5-2\rho_6\right)\tilde{J}_2-\left(2\rho_3+\rho_4+2\rho_1\right)\tilde{J}_3\\
27M^{33} &= \left(\tilde{J}_1-\tilde{J}_3+\tilde{J}_2\right)4\rho_3+\left(2\rho_8+4\rho_7+2\rho_6+\rho_5+4\rho_4+8\rho_3+4\rho_2+2\rho_1\right)J_0+\\
&\quad +\left(\rho_6+2\rho_7+2\rho_2\right)\tilde{J}_1+\left(\rho_8+2\rho_4+2\rho_7\right)\tilde{J}_2-\left(2\rho_2+2\rho_4+\rho_1\right)\tilde{J}_3\\
27M^{44} &= \left(-\tilde{J}_1-\tilde{J}_3+\tilde{J}_2\right)4\rho_4+\left(4\rho_8+2\rho_7+\rho_6+2\rho_5+8\rho_4+4\rho_3+2\rho_2+4\rho_1\right)J_0+\\
&\quad +\left(-\rho_5-2\rho_8-2\rho_1\right)\tilde{J}_1+\left(\rho_7+2\rho_3+2\rho_8\right)\tilde{J}_2-\left(\rho_2+2\rho_3+2\rho_1\right)\tilde{J}_3
\end{aligned}$$

$$
\begin{aligned}
27M^{55} &= \left(-\tilde{J}_1 - \tilde{J}_2 + \tilde{J}_3\right)4\rho_5 + \left(4\rho_8 + 2\rho_7 + 4\rho_6 + 8\rho_5 + 2\rho_4 + \rho_3 + 2\rho_2 + 4\rho_1\right)J_0 + \\
&\qquad + \left(-\rho_4 - 2\rho_8 - 2\rho_1\right)\tilde{J}_1 + \left(-\rho_2 - 2\rho_1 - 2\rho_6\right)\tilde{J}_2 + \left(2\rho_8 + \rho_7 + 2\rho_6\right)\tilde{J}_3 \\
27M^{66} &= \left(\tilde{J}_1 - \tilde{J}_2 + \tilde{J}_3\right)4\rho_6 + \left(2\rho_8 + 4\rho_7 + 8\rho_6 + 4\rho_5 + \rho_4 + 2\rho_3 + 4\rho_2 + 2\rho_1\right)J_0 + \\
&\qquad + \left(\rho_3 + 2\rho_7 + 2\rho_2\right)\tilde{J}_1 + \left(-2\rho_2 - 2\rho_5 - \rho_1\right)\tilde{J}_2 + \left(2\rho_5 + 2\rho_7 + \rho_8\right)\tilde{J}_3 \\
27M^{77} &= \left(J_2 + J_1 + J_3\right)4\rho_7 + \left(4\rho_8 + 8\rho_7 + 4\rho_6 + 2\rho_5 + 2\rho_4 + 4\rho_3 + 2\rho_2 + \rho_1\right)J_0 + \\
&\qquad + \left(2\rho_6 + \rho_2 + 2\rho_3\right)\tilde{J}_1 + \left(2\rho_3 + \rho_4 + 2\rho_8\right)\tilde{J}_2 + \left(2\rho_6 + \rho_5 + 2\rho_8\right)\tilde{J}_3 \\
27M^{88} &= \left(\tilde{J}_2 - \tilde{J}_1 + \tilde{J}_3\right)4\rho_8 + \left(8\rho_8 + 4\rho_7 + 2\rho_6 + 4\rho_5 + 4\rho_4 + 2\rho_3 + \rho_2 + 2\rho_1\right)J_0 + \\
&\qquad + \left(-\rho_1 - 2\rho_4 - 2\rho_5\right)\tilde{J}_1 + \left(2\rho_7 + \rho_3 + 2\rho_4\right)\tilde{J}_2 + \left(2\rho_5 + \rho_6 + 2\rho_7\right)\tilde{J}_3
\end{aligned}
\tag{20}
$$

CM semi-analytical closed-form integration rule is summarized by

| |
|---|
| compute $J^0_{mn}$ given by (12) |
| compute $J_0$ using (5) |
| compute $M^{ii}$ given by (13) |

LM semi-analytical closed-form integration rule is summarized by

| |
|---|
| compute $J^k_{mn}\,(k=0,1,2,3)$ given by (12)(17)(18)(19) |
| compute $J_0, J_{point1}, J_{point2}, J_{point3}$ using (5) |
| compute $\tilde{J}_k$ given by (15) |
| compute $M^{ii}$ use (20) |

Generally speaking, generation of integration rule using the proposed approach contains two steps: The first is to develop a model for jacobian determinant namely $J = \tilde{N}^i(\xi,\eta,\zeta)\tilde{J}_i(N_{jk})\ (i = 0,..,n_J)$ where $\tilde{N}^i$ are coordinate dependent ansatz functions, e.g. polynomial or monomial terms, and $\tilde{J}_i$ nodal component dependent terms. The second step is analytical integration of the mass matrix components.

## 4. Preliminary numerical study.

In no way the present letter pretends to have a complete, all inclusive, numerical study; however the preliminary numerical study illuminates obvious benefits of using CM and LM semi-analytical rules over equivalent schemes.

Specific values of initial nodal densities are given by

$$\rho_1 = 1 \ , \ \rho_2 = 1 \ , \ \rho_3 = 1 \ , \ \rho_4 = 1 \ , \ \rho_5 = 2 \ , \ \rho_6 = 2 \ , \ \rho_7 = 2 \ , \ \rho_8 = 2 \tag{21}$$

Consider the next parallelepiped element

$$\begin{aligned}
&N_{1,1} = -1+\varepsilon \, , N_{1,2} = -1 \, , N_{1,3} = -1 \, , N_{2,1} = 1+\varepsilon \, , N_{2,2} = -1 \, , N_{2,3} = -1 \\
&N_{3,1} = 1 \, , N_{3,2} = 1 \, , N_{3,3} = -1+\varepsilon \, , N_{4,1} = -1 \, , N_{4,2} = 1 \, , N_{4,3} = -1+\varepsilon \\
&N_{5,1} = -1+\varepsilon \, , N_{5,2} = -1 \, , N_{5,3} = 1 \, , N_{6,1} = 1+\varepsilon \, , N_{6,2} = -1 \, , N_{6,3} = 1 \\
&N_{7,1} = 1 \, , N_{7,2} = 1 \, , N_{7,3} = 1+\varepsilon \, , N_{8,1} = -1 \, , N_{8,2} = 1 \, , N_{8,3} = 1+\varepsilon
\end{aligned} \tag{22}$$

For $\varepsilon = 0$ the above yield a cube with edge length equal 2, although for $\varepsilon = 100$, (22) lead to a very skewed non-rectangular parallelepiped element with big aspect ratio. Lumped mass matrix components are calculated and absolute error is estimated with respect to exact values, then absolute error is averaged between 8 components.

| $\varepsilon = 100$ | CM | $n_p = 1$ | LM | $n_p = 4$ |
|---|---|---|---|---|
| **Error %** | 0 | 11.25 | 0 | 0 |

For non-homogeneous non-rectangular parallelepiped element CM over perform its equivalent, numerical integration based on one point quadrature (9)(10).

We want to examine accuracy performance for a coarse mesh. Consider the next element family

$$\begin{aligned}
&N_{1,1} = -1+R \, , N_{1,2} = -1+R \, , N_{1,3} = -1+R \, , N_{2,1} = 1+R \, , N_{2,2} = -1+R \, , N_{2,3} = -1+R \\
&N_{3,1} = 1+R \, , N_{3,2} = 1+R \, , N_{3,3} = -1+R \, , N_{4,1} = -1+R \, , N_{4,2} = 1+R \, , N_{4,3} = -1+R \\
&N_{5,1} = -1+R \, , N_{5,2} = -1+R \, , N_{5,3} = 1+R \, , N_{6,1} = 1+R \, , N_{6,2} = -1+R \, , N_{6,3} = 1+R \\
&N_{7,1} = 1+R \, , N_{7,2} = 1+R \, , N_{7,3} = 1+R \, , N_{8,1} = -1+R \, , N_{8,2} = 1+R \, , N_{8,3} = 1+R
\end{aligned} \tag{23}$$

Where $R$ is a random variable which is uniformly distributed between $-\delta$ and $\delta$. Pseudo-random numbers are produced with MAPLE™ built in function. For each component of every element $R$ takes different real number in the range. For $\delta = 0$ element family (23) reduces to a cube with edge length equal 2. We've studied $\delta$ values in the range $0 \le \delta \le 0.7$. For each value of $\delta$ one hundred different elements has been produced. For every element an exact and approximate lumped mass matrix components were computed. Averaged absolute error results are reported in Figure 1.

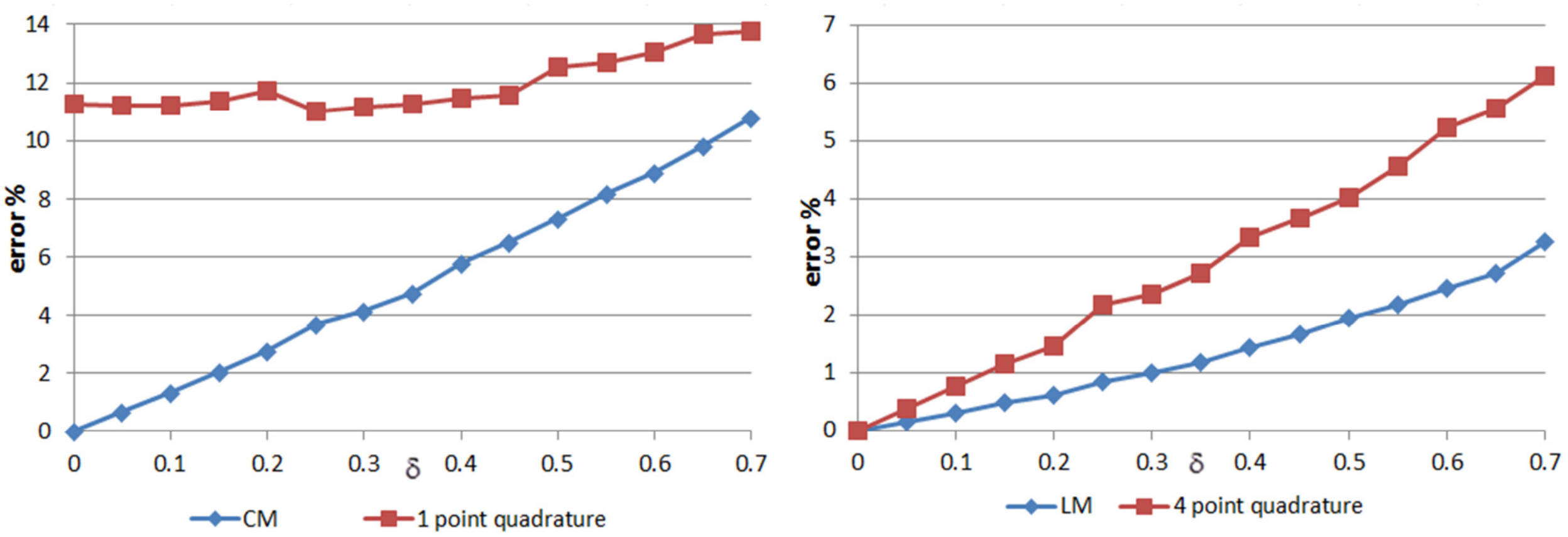

Figure 1: Averaged absolute error is presented as a function of delta. For each point 100 elements are used. Left graph is showing the CM rule and one point quadrature numerical integration rule. Right graph is showing the LM rule vs numerical integration with four point quadrature.

## 5. Conclusions

In this study, for the first time, two low order semi-analytical integration rules for mass matrix of an 8-node brick element are discussed. CM assumes constant jacobian determinant while LM allows linear variation of the metric in element domain. Mass matrix component follow from analytical integration.

Both closed-form schemes are exact for parallelepiped mesh. Preliminary numerical study for coarse mesh is conducted. Random mesh is generated such that one parameter delta controls the coarseness of the mesh. Preliminary numerical study established that the averaged absolute error is always lower than for equivalent schemes.